\documentclass[11pt]{article}

\usepackage{a4wide}
\def\versiondate{4. Dez.~2005}
\usepackage{amssymb}
\usepackage{amsmath}

\usepackage{graphicx}
\usepackage{color}
\graphicspath{{Figures/}}
\input{psfig.sty}

\def\comment#1{}

\def\NI{\noindent}
\def\ni{\noindent}
 
\def\sk{\smallskip}

\def\term#1{{\em #1}}

\def\ITEMMACRO #1 ??? #2 ???{\par\vskip4pt\noindent%
\hangindent=#2em\setbox0\hbox{#1\kern4pt}%
\ifdim\wd0<\hangindent\setbox0\hbox to\hangindent{\hss#1\kern7pt}\fi%
\box0\ignorespaces}
 
\def\Item(#1){\ITEMMACRO {\rm (#1)} ??? 1.8 ???}
\def\FreeItem#1{\ITEMMACRO {#1} ??? 1.8 ???}

\def\bItem{\ITEMMACRO \hss$\bullet$ ??? 1.8 ???}
\let\Bitem=\bItem
\def\BrackItem[#1]{\ITEMMACRO [#1] ??? 1.8 ???}

\newtheorem{theorem}{Theorem}
\newtheorem{lemma}{Lemma}

\newtheorem{proposition}{Proposition}

\newtheorem{corollary}{Corollary}

\def\Fact#1.{\par\sk{\NI\bf Fact~#1}\ }
\def\Problem#1.{\par\sk{\NI\bf Problem~#1}\ }
\def\Claim#1.{\medbreak\ni{\bf Claim~#1}\ }
\def\Case#1.{\medbreak\ni{\bf Case~#1}\ }
\def\SubCase#1.{\medbreak\ni{\bf Subcase~#1}\ }

\def\Proof{\ni{\sl Proof}\ }
\def\qed{\hfill\fbox{\hbox{}}\medskip}
\def\qedclaim{\hfill$\triangle$\smallskip}

\def\below#1{[\kern1pt#1-1\kern-1.5pt\updownarrow\kern-1.5pt#1\kern1pt]}
\def\between#1#2{[\kern1pt#1\kern-1.5pt\updownarrow\kern-1.5pt#2\kern1pt]}
\def\bbelow#1{[\kern1pt#1-1\kern-1.5pt\updownarrow\kern-1pt\updownarrow\kern-1.5pt#1\kern1pt]}
\def\bbetween#1#2{[\kern1pt#1\kern-1.5pt\updownarrow\kern-1pt\updownarrow\kern-1.5pt#2\kern1pt]}

\def\ol{\overline}
\def\with{\kern1.5pt :\kern1.5pt}
\def\spanned{\hbox{\sf span}_2}
\def\span#1{\hbox{\sf span}_#1}
\def\weaksp{\hbox{\sf span}_4}
\def\exposed{\hbox{\sf exposed}}
\def\RN{R_{\bf 0}}

\def\symbolboxwidth{\hbox to 3.2pt}
\def\symbolrulewidth{\hrule width3.25pt}
\def\symbolboxheight{\vbox to 3pt}

\def\nw{{\vbox{\symbolboxwidth{\hfil\symbolboxheight{\vfil}\vrule}\symbolrulewidth}}}
\def\se{{\vbox{\symbolrulewidth\symbolboxwidth{\vrule\hfil\symbolboxheight{\vfil}}}}}
\def\ne{{\vbox{\symbolboxwidth{\vrule\hfil\symbolboxheight{\vfil}}\hrule }}}
\def\sw{{\vbox{\symbolrulewidth\symbolboxwidth{\hfil\symbolboxheight{\vfil}\vrule}}}}
\def\rec{{\vbox{\symbolrulewidth\symbolboxwidth
         {\vrule\hfil\symbolboxheight{\vfil}\vrule}\symbolrulewidth}}}

\def\join{{\vee}}

\def\<{\langle}
\def\>{\rangle}

\def\RR{\hbox{\sf I\kern-1ptR}}
\def\NN{\hbox{\sf I\kern-1ptN}}
\def\ZZ{\hbox{\sf Z\kern-4ptZ}}

\begin{document}
%

\begingroup\raggedright
\baselineskip22pt
\LARGE\bf

Empty Rectangles and Graph Dimension

\endgroup
\vskip20pt

\begingroup
\fontsize{12}{14}%
{\sc Stefan Felsner}
\vskip6pt

{\small\it Technische Universit\"at Berlin,  Institut f\"ur Mathematik, MA 6-1,}
\hfill\vskip-1pt
{\small\it Stra\ss{}e des 17. Juni 136, 10623 Berlin, Germany}\hfill\vskip1pt
{\small {\it E-mail:} felsner@math.tu-berlin.de}
\endgroup

\footnotetext{\versiondate}


\vskip20pt
\begingroup\fontsize{10}{12}\rm
\noindent{\bf Abstract}
We consider rectangle graphs whose edges are defined by pairs of points in
diagonally opposite corners of empty axis-aligned rectangles.
The maximum number of edges
of such a graph on $n$ points is shown to be $\lfloor\frac{1}{4}n^2 +n -2\rfloor$.
This number also has other interpretations:
\bItem
It is the maximum number of edges of a graph of dimension $\bbetween{3}{4}$,
i.e., of a graph with a realizer of the form $\pi_1,\pi_2,\ol{\pi_1},\ol{\pi_2}$.
\bItem
It is the number of 1-faces in a special Scarf complex.
\par
\medskip

\ni
The last of these interpretations allows to deduce the maximum number of 
empty axis-aligned rectangles spanned by 4-element subsets of a set of $n$ points.
Moreover, it follows that the extremal point sets for the two problems coincide.

We investigate the maximum number of of edges of a graph of
dimension $\between{3}{4}$, i.e., of a graph with a realizer of the form
$\pi_1,\pi_2,\pi_3,\ol{\pi_3}$. 
This maximum is shown to be $\frac{1}{4}n^2 + O(n)$.

Box graphs are defined as the 3-dimensional analog of rectangle
graphs. The maximum number of edges of such
a graph on $n$ points is shown to be $\frac{7}{16}n^2  + o(n^2)$.

\vskip10pt\ni
\endgroup

\section{Introduction}         

A set of points in the plane can serve as the vertex set of various
geometrically defined graphs. The most popular example is the Delaunay
triangulation which has an edge between two points iff there is an
empty circle through them. Graph theoretically Delaunay triangulations
are a subclass of planar triangulations. Schnyder showed that the set of all planar
triangulations can be obtained if we let two points form an edge iff
they are on the boundary of an empty triangle with sides parallel to
three given lines.

In this paper we consider rectangle graphs whose edges are defined by
pairs of points on the boundary of empty axis-aligned rectangles and
variants of it.  An appealing aspect of the topic is that the objects
dealt with have natural interpretations in different areas. The first
three sections approach the theme from different directions and
establish connections between them. These sections are
\FreeItem{2.} Empty Rectangles and Empty Boxes
\FreeItem{3.} Dimension of Graphs
\FreeItem{4.} Orthogonal Surfaces and Scarf's Theorem
\medskip

\ni
The lengthy proofs of
Theorems~\ref{thm:span} and~\ref{thm:singlearrow} are taken to the
appendix.

\section{Empty Rectangles and Empty Boxes}

Throughout this paper a \term{rectangle} is an axis-aligned
rectangle in the plane $\RR^2$, i.e., a set of the form
$$
R = R(r_1,r'_1,r_2,r_2') = \big\{ s=(s_1,s_2) \in \RR^2 \with r_1 \leq s_1
\leq r'_1 \hbox{ and } r_2 \leq s_2 \leq r'_2 \big\}.
$$
Given a finite set $X$ of points
a rectangle $R$ is an \term{empty rectangle} if the open interior
$R^o$ of $R$ contains no point of $X$. 
The rectangle $R$ \term{spanned} by $A \subseteq X$ is the smallest
rectangle containing all points of $A$,
\begin{eqnarray*}
& R = R[A]  =  \big\{ s \in \RR^2 \with & 
         \min(x_1 : x \in A) \leq s_1 \leq \max(x_1 : x\in A) \hbox{ and }\\
         & & \min(x_2 : x\in A) \leq s_2 \leq \max(x_2 : x \in A) \big\}.
\end{eqnarray*}
Define the \term{rectangle graph $G_r(X) =(X,E_r)$ of the point set $X$} such
that two points $p,q\in X$ are an edge in $E_r$ iff they span an empty
rectangle $R[p,q]$. The edges of the rectangle graph $G_r(X)$
are pairs of points from $X$ spanning an empty rectangle.
Let $\spanned(X)$ be the number of edges of $G_r(X)$.
In the degenerate situation where all points of $X$ have one
coordinate in common the graph $G_r(X)$ is a complete graph.
A set $X$ is called \term{generic} if no two
points of $X$ share a coordinate. 

\begin{theorem}\label{thm:span}
For every generic set $X$ of $n$ points in $\RR^2$
$$
\spanned(X) \leq \Big\lfloor\frac{n^2}{4} +n -2\Big\rfloor,
$$
this bound is best possible.
\end{theorem}

\ni
The proof is in Section~\ref{sec:proof1}.
\medskip

Generalizing notation we let $\span{k}(X)$ be the number of
$k$-subsets $A$ of $X$ such that $R[A]$ is an empty rectangle, that is
the interior $R^o$ is empty and all elements of $A$ are on the
boundary of $R$.  If $X$ is generic, then each side of a rectangle can
contain at most one point, hence $\span{k}(X) = 0$ for $k\geq 5$. The
numbers $\span{3}(X)$ and $\span{4}(X)$ count the triangles and
4-cliques of $G_r(X)$.

A point $p$ in $X$ is \term{orthogonally exposed} if
one of the four quadrants determined by $p$ is empty, i.e., contains
no points from $X$.  Let $\exposed(X)$ be the number of orthogonally
exposed points of $X$. There are linear relations between the
quantities $\spanned(X)$, $\span{4}(X)$ and $\exposed(X)$
as well as $\spanned(X)$, $\span{3}(X)$ and $\exposed(X)$.

\begin{theorem}\label{thm:w-span}
For every generic set $X$  in $\RR^2$
\begin{eqnarray*}
\spanned(X) - \span{4}(X) + \exposed(X) &=&  3\,(|X|-1),\\
2\;\spanned(X) - \span{3}(X) + \exposed(X) &=&  4\,(|X|-1).
\end{eqnarray*}
\end{theorem}

\PsFigCap{48}{example-small}
{A point set $X$ with its graph $G_r(X)$. The statistics are
$|X|=6$, $\exposed(X)=5$, $\spanned(X)=12$, $\span{3}(X)=9$ and $\span{4}(X)=2$.}

I include the following proposition because I think that
the result is surprising, at least at the first glance.
The proof is left to the reader as an exercise.

\begin{proposition}
Let $X$ be a set of $n$ random points from the unit square.
The expected number of edges of the rectangle graph $G_r(X)$
equals the expected number of comparisons performed
by random quick-sort on an $n$ element input, i.e., it is
$\sum_{1\leq i<j\leq n}\frac{2}{j-i+1}$. 
\end{proposition}

\subsection{Diagrams of 2-Dimensional Orders}

An order $P=(X,<_P)$ is 2-dimensional if it can be represented as
dominance order on a generic set of points in $\RR^2$, i.e.,
$(x_1,x_2) <_P (y_1,y_2)$ iff $x_1 < y_1$ and $x_2 < y_2$.
The edges in the diagram of $P$ then correspond to spanned
empty rectangles with one point in the lower-left and one point in the
upper-right corner. 

A conjugate $P^c$ of $P$ can be obtained from the representation
by reverting the direction of one axis. If the first coordinate was
reverted, then we have  $(x_1,x_2) <_{P^c} (y_1,y_2)$ in $P^c$ iff
$x_1 > y_1$ and $x_2 < y_2$. The edges in the diagram of $P^c$ correspond to spanned
empty rectangles with one point in the upper-left and one point in the
lower-right corner. Hence, the union of $P$ and $P^c$ is the
rectangle graph $G_r(X)$

Since diagrams are triangle free graphs, they have at most $n^2/4$
edges, this bound is tight even for 2-dimensional
orders. Theorem~\ref{thm:span} shows that the diagrams of $P$ and a
conjugate $P^c$ together never have much more edges.

The order theoretic point of view can yield insights about 
rectangle graphs. Kr{\'i}z and Ne\v{s}et\v{r}il \cite{kn-cnhded-91} 
have shown that diagrams of 2-dimensional orders have unbounded
chromatic number, hence, so do rectangle graphs. 
R\"odel and Winkler showed that almost all 2-dimensional orders enjoy
this property since their diagrams have independence number
$o(n)$. The proof is lost and Peter Winkler says:
``So we are obliged to tell people that the problem is again open!''

\subsection{Box graphs}

Define box graphs as 3-dimensional analogs to rectangle graphs.
That is, for a set $Y$ of points in $\RR^3$ we define the
{\it box graph $G_b(Y)$} such
that two points $p,q\in Y$ form an edge iff they span an empty
axis aligned box $B[p,q]$. We ask the question for the
maximum number of edges of a box graph defined by $n$ points.

\begin{theorem}
Box graphs of sets of $n$ point have at most $\frac{7}{8}n^2 +o(n^2)$
edges. Up to the error term this bound is best possible.
\end{theorem}

\Proof
In the next section we prove:
\Item(1)
$K_{17}$ is a forbidden subgraph for the class of box graphs
(Proposition~\ref{prop:K17}).
\Item(2)
Every complete 8-partite graph is a subgraph of some box graph
(Proposition~\ref{prop:8part}). 
\Item(3)
$K_{16}$ is a box graph (Proposition~\ref{prop:K16}).
\Item(4)
There are complete 9-partite graphs which are forbidden as
subgraphs of box graphs.
(Proposition~\ref{prop:no9part}). 
\medskip

\ni
With (2) we have box graphs with at least $\frac{7}{8}n^2$ edges.
For the upper bound we need the Erd\H{o}s-Stone Theorem which asserts
that for any fixed graph $H$ with chromatic number $\chi(H)=k$ the
number of edges of graphs containing no subgraph isomorphic to $H$ is
of order $\frac{k-2}{2(k-1)}n^2 + o(n^2)$. Since the forbidden
subgraphs from (4) have chromatic number 9 we obtain the 
bound.
\qed

\section{Dimension of Graphs}

Let $G = (V,E)$ be a finite simple graph.
A nonempty family $\mathcal{R}$ of permutations (linear orders) of the
vertices of a graph $G$ is called a \term{realizer} 
of $G$ provided

\Item($*$)  For every edge $S\in E$ and every vertex $x\in V\setminus S$,
there is some $\pi\in\mathcal{R}$ so that $x>y$ in $\pi$ for
every $y\in S$.
\par

\medskip\ni
The \term{dimension} of $G$, denoted $\dim(G)$, is defined 
as the least positive integer $t$ for which $G$ has a realizer of 
cardinality $t$.  
\medskip


In order to avoid trivial complications when the condition
($*$) is vacuous, we restrict our attention to connected graphs with
three or more vertices.

The definition of dimension has a natural generalization
to hypergraphs: Just replace `edge' by `hyperedge' in condition
($*$). We review two, by now classical, facts about the dimension of graphs
and hypergraphs:

\bItem A graph $G$ is planar if and only if
its dimension is at most~$3$.\\
(Schnyder~\cite{s-pgpd-89})

\bItem The hypergraph $H=(V,F)$ of vertices $V$ versus edges and faces 
of a 3-polytope has dimension 4 but the dimension drops to 3 
when a face $f$ is removed from $F$.\\
(Brightwell-Trotter~\cite{bt-odcp-93} --
see \cite{f-cdpgo-01}, \cite{f-gepg-03} and \cite{fz-os3d-05} for simplified proofs)
\medskip

\ni In~\cite{ft-ppg-05} a refined concept for dimension of graphs was
investigated.  For an integer $t\ge2$, define the dimension of a graph
to be $\between{t}{t+1}$ if it has dimension greater than $t$ yet has
a realizer of the form $\{\pi_1,\pi_2,\dots,\pi_t,\pi_{t+1}\}$ with
$\pi_{t+1}=\ol{\pi_{t}}$, where $\ol{\pi}$ is the \term{reverse} of $\pi$,
i.e., $x<y$ in $\pi$ if and only if $x>y$ in $\ol{\pi}$ for all $x,y\in
X$.  The following facts about this `intermediate' dimension are known
(Felsner-Trotter~\cite{ft-ppg-05}):

\bItem A graph $G$ is outerplanar if and only if
its dimension is at most~$\between{2}{3}$.

\bItem A graph $G$ with dimension at most $\between{3}{4}$
has at most $\frac{1}{4}n^2 + o(n^2)$ edges.
\medskip

\ni
Generalizing the notation we say that the dimension of a graph 
is at most $\bbetween{t-1}{t}$ if it has a realizer of the form 
$\{\pi_1,\pi_2,\dots,\pi_t\}$ with $\pi_{t}=\ol{\pi_{t-2}}$
and $\pi_{t-1}=\ol{\pi_{t-3}}$. 

\begin{proposition}\label{prop:doublearrow}
Graphs of dimension at most $\bbetween{3}{4}$
are exactly the rectangle graphs.
\end{proposition}

\Proof
We just have to translate the edge-condition ($*$) for a realizer
$\{\pi_1,\pi_2,\ol{\pi_1},\ol{\pi_2}\}$ into the emptiness of a spanned
rectangle and vice versa.

Given a set $X$ of points in the plane $\RR^2$ we obtain two linear
orders $\pi_1(X)$ and $\pi_2(X)$ from projections to the coordinate axes.
Conversely, the first two linear orders of a realizer determine a pair
of integer coordinates in the range $[1,\ldots,n]$ and thus a set $X$
of points with $|X|=n$.

Let $p,q\in X$ span an empty rectangle.  This implies
that for every $z \in X\setminus\{p,q\}$ at least one of the following
four relations must be false: $(1)$~$z_1 < \max(p_1,q_1)$ or
$(2)$~$z_2 < \max(p_2,q_2)$ or $(3)$~$z_1 > \min(p_1,q_1)$
or $(4)$~$z_2 > \min(p_2,q_2)$. If the $i$th of these relations is false, then
$z > p,q$ in $\pi_i$, where $\pi_{3}=\ol{\pi_{1}}$
and $\pi_{4}=\ol{\pi_{2}}$. This in  turn implies condition ($*$) for the 
edge $\{p,q\}$. 
Since the argument can be reversed we obtain: 
$p,q$ span an empty rectangle exactly if
$p,q$ satisfy condition ($*$).
\qed

From the proposition and Theorem~\ref{thm:span} 
we obtain:

\begin{corollary}
A graph $G$ of dimension at most $\bbetween{3}{4}$
has at most $\lfloor\frac{1}{4}n^2 + n -2 \rfloor$ edges.
\end{corollary}

Geometric arguments enable us to improve upon the upper bound of
order $\frac{1}{4}n^2 + o(n^2)$ (see~\cite{ft-ppg-05}) for the number
of edges of a graph of dimension $\between{3}{4}$.

\begin{theorem}\label{thm:singlearrow}
A graph $G$ of dimension at most $\between{3}{4}$
has at most $\frac{1}{4}n^2 + O(n)$ edges.
\end{theorem}

From the proof of the theorem in Section~\ref{sec:proof2} it follows
that actually $\frac{1}{4}n^2 + 5n$ is a valid bound for all~$n$.
This is very much in contrast to the old bound which was obtained by
combining the Product Ramsey Theorem and the Erd\H{o}s-Stone Theorem.
Therefore, the resulting bound was only asymptotic.

In the preprint preceeding the publication~\cite{ft-ppg-05}
a structural characterization of graphs with dimension at most
$\between{3}{4}$ was conjectured. This conjecture would have
implied that these graphs have at most $\frac{1}{4}n^2 + 2n - 6$
edges, a construction with that number of edges is known.
Though the conjecture was disproved in~\cite{odmr-hd-05} it is still
possible that the bound on the number of edges is correct.

\Problem 1.{\it
Is it true that graphs of dimension at most $\between{3}{4}$
have at most $\frac{1}{4}n^2 + 2n - 6$ edges?
}

In \cite{aft-mnegb-99} it has been shown that a graph of dimension
$\leq 4$ (i.e., a graph admitting a realizer
$\{\pi_1,\pi_2,\pi_3,\pi_4\}$ ) has at most $\frac{3}{8}n^2 + o(n^2)$
edges.  The proof combines the Product Ramsey Theorem and the
Erd\H{o}s-Stone Theorem. We will use this technique in the next
subsection when we bound the number of edges of box graphs.
 
It would be interesting to have a proof of the $\frac{3}{8}n^2$ bound which
is elementary and geometric. Such a proof should also yield
an improvement in the error term.

\Problem 2. Show that a graph with a realizer
$\{\pi_1,\pi_2,\pi_3,\pi_4\}$ has at most $\frac{3}{8}n^2 + Cn$
edges for some reasonable constant $C$.
\medskip

\subsection{Realizer for box graphs}\label{sec:3-6}

We now come back to box graphs. Just as in the proof of
Proposition~\ref{prop:doublearrow} it can be shown that
box graphs are exactly the graphs with a realizer
$\{\pi_1,\pi_2,\pi_3,\ol{\pi_1},\ol{\pi_2},\ol{\pi_3}\}$.

\begin{proposition}\label{prop:K17}
Box graphs contain no 17-cliques.
\end{proposition}

\Proof
Recall the Erd\H{o}s-Szekeres Theorem: A permutation of $n$ numbers
contains a monotone subsequence of length $\lceil\sqrt{n}\rceil$.

Let $G$ be a graph with vertices labeled $1$ to $n$ and a realizer
$\{\pi_1,\pi_2,\pi_3,\ol{\pi_1},\ol{\pi_2},\ol{\pi_3}\}$.
We may assume that $\pi_1$ is the identity permutation.
Consider a set $A$ of 17 vertices of $G$. Permutation
$\pi_2$ contains a monotone subsequence of length at least 5
on the vertices from $A$, let $B$ be such a set of 5 vertices.
Permutation
$\pi_3$ contains a monotone subsequence of length at least 3
on the vertices from $B$, let $\{i,j,k\}$ be such a set of 3 vertices.
The vertices $i,j,k$ are in this order or in the reverse in each of
the six permutations of the realizer. Hence, $j$ obstructs the edge
$i,k$ and the set $A$ can't induce a clique.
\qed
%
%
\begin{proposition}\label{prop:8part}
Every complete 8-partite graph is a subgraph of some box graph.
\end{proposition}

\Proof
Consider pairwise disjoint sets $X_i$ for $i=1,..,8$ and let
$\sigma_i$ be a permutation of $X_i$. Define permutations
by concatenations of the $\sigma_i$'s as follows:
\begin{eqnarray*}
\pi_1 & = &
\sigma_1+\sigma_2+\sigma_3+\sigma_4+\sigma_5+\sigma_6+\sigma_7+\sigma_8\\
 \pi_2 & = & 
\sigma_5+\sigma_3+\sigma_2+\sigma_8+\sigma_1+\sigma_7+\sigma_6+\sigma_4\\
\pi_3 & = & 
\sigma_7+\sigma_4+\sigma_8+\sigma_6+\sigma_3+\sigma_1+\sigma_5+\sigma_2
\end{eqnarray*}
The claim is that these permutations together with the reversed form a
realizer of a graph containing all edges between vertices from
different sets $X_i$. Since each permutation comes with its reversed
the condition $(*)$ for edges can be written as

\Item($**$)  For every $x,y$ with $x\in X_i$, $y\in X_j$, $i\neq j$
and  every vertex $z \neq x,y$, there is an $\pi_a$ so that
$z$ is not between $x$ and $y$ in $\pi_a$.
\smallskip

\ni
If $z\in X_k$ with $k\neq i,j$, then $(**)$ is satisfied if
$\sigma_i$ and $\sigma_j$ are consecutive in one of the permutations.
This condition is fulfilled for all pairs $\{i,j\}$ except
$1,4$ and $1,6$ and $2,4$ and $2,6$ and $2,7$ and $3,7$ and $3,8$
and $5,7$ and $5,8$. For these pairs it is enough to look at $\pi_1$
and $\pi_2$ to verify that there is no $\sigma_k$ which is between
$\sigma_i$ and $\sigma_j$ in all three permutations.

Now consider $z\in X_i$. If we have $\sigma_i$ before $\sigma_j$ in
$\pi_a$ and $\sigma_j$ before $\sigma_i$ in
$\pi_b$, then $z$ is outside the interval $x,y$ in either $\pi_a$ or
$\pi_b$. Again it is not hard to inspect all pairs $i,j$ of indices
and verify this criterion for them.

\begin{proposition}\label{prop:K16}
$K_{16}$ is a box graph.
\end{proposition}

\Proof
Do the above construction from the proof of the previous
proposition with sets $X_i$ of
cardinality 2 each. The pair of vertices from $X_i$ stays together in 
all the permutations, hence, the edge condition $(**)$ is trivially
satisfied for these pairs.
\qed
 
\begin{proposition}\label{prop:no9part}
For sufficiently large $n$ the complete  9-partite graph
${\bf T}(9n,9)$  with parts
of equal size $n$ is a forbiden subgraph for box graphs.
\end{proposition}

For the proof we need a Ramsey type theorem. This Theorem
has been used in \cite{aft-mnegb-99} and \cite{ft-ppg-05} in a similar
context.
\smallskip

\ni
{\bf Product Ramsey Theorem\ }\bgroup\it 
Given positive integers $m$, $r$ and $t$, there exists
an integer $n_0$ so that if 
$S_1,\dots,S_t$ are sets with $|S_i|\ge n_0$ for
all $i$, and $f$ is any map
which assigns to each transversal $(x_1,\ldots,x_t)$ with $x_i\in S_i$
 a color from $[r]$,
then there exist $H_1,\dots,H_t$ with $H_i\subseteq S_i$ for all $i$ and a color
$\alpha\in[r]$ so that
\bItem $f(g)=\alpha$ for every transversal $g=(x_1,\ldots,x_t)$ 
with $x_i\in H_i$, for all $i=1,2,\dots,t$.
\egroup
\bigskip

\ni
{\sl Proof of Proposition~\ref{prop:no9part}.\ }
Let ${\bf T} = {\bf T}(9n,9)$ be the complete 9-partite graph with parts
of equal size $n$. Suppose that there is a box graph $G$ containing
${\bf T}$ as a subgraph.
Let $S_1,\ldots,S_9$ be the independent sets  of
${\bf T}$. Define a coloring of the transversals as follows.  
Each transversal is just a~$9$-element set containing
one point from each $S_i$.
Then consider the order of these~$9$ points
in the permutations $\pi_1$, $\pi_2$, $\pi_3$ of a realizer for $G$.
In each $\pi_a$, the
$9$ points can occur in any of $9!$ orders.  So taking the
$3$ orders altogether, there are at $r=(9!)^3$ patterns.

Applying the Product Ramsey Theorem, it follows that, if $n$ was large
enough, then for each $i\in[9]$, there is a $2$-element subset $H_i$
of $S_i$ so that all the transversals with elements from these subsets
receive the same color.  This implies that the linear orders treat the
sets $H_1$, $H_2$, $H_3$, $H_4$ and $H_5$ as \textit{blocks}, i.e., if
a point from one block is over a point from another block in $\pi_k$,
then both points from the first block are over both points from the
second block in $\pi_k$.

Now restrict the realizer to the vertex set $H_1\cup H_2\cup \dots \cup
H_9$. This gives a graph with all the edges of ${\bf T}(18,9)$ but
since the vertices of each $H_i$ stay together in 
all the permutations they also satisfy the edge condition $(**)$.
Therfore, the restricted realizer represents
a complete graph $K_{18}$. This is a contradiction since 
 $K_{18}$ is not a box graph (Proposition~\ref{prop:K17}).
\qed

\section{Orthogonal Surfaces and Scarf's Theorem}

Let $V\subset \RR^d$ be a (finite) antichain $V$ in the dominance order
on $\RR^d$. The \term{orthogonal surface} $S_V$ generated by $V$ is
the topological boundary of the filter  $\< V \> = \big\{ x \in \RR^d
: \exists v\in V \text{ such that } x_i \geq v_i \text{ for } i=1,..,n
\big\}$.

An orthogonal surface $S_V$ in $\RR^d$ is \term{suspended} 
if $V$ contains one element from each positive coordinate axis and the
coordinates of all the other members of $V$ are strictly positive.
An orthogonal surface $S_V$ is \term{generic} if no two points in $V$
have a common coordinate, i.e., $v_i \neq v'_i$ for all $v,v'\in V$
and $i=1,..,d$. In the case of a suspended surface the genericity
condition is only applied to coordinates of positive value.

The \term{Scarf complex} $\Delta_V$ of a generic orthogonal surface $S_V$
generated by $V$ consists of all the subsets $U$ of $V$ with the
property that $\bigvee\{ v : v\in U \} \in S_V$, the join $\bigvee\{ v : v\in U \}$
is the point $u$ with coordinates $u_i = \max(v_i : v\in U)$.
The property $u\in S_V$ is equivalent to either of $i$ and $ii$:

\Item($i$) There is no $v\in V$ which is strictly dominated by $u$
(a point $q$ is strictly dominated by $p$ if $q_i < p_i$ for
$i=1,..,d$).
\par

\Item($ii$) Every $v\in V$ which is dominated by $u$
has at least one coordinate in common with $u$.
\par

\medskip\ni
It is a good exercise to show that the Scarf complex $\Delta_V$ 
of a generic $V$ is a simplicial complex.

\begin{theorem}[Scarf'73]\label{thm:scarf}
The Scarf complex $\Delta_V$ of a generic suspended orthogonal surface $S_V$
in $\RR^d$ is isomorphic to the face complex of a simplicial
$d$-polytope with one facet removed.
\end{theorem}

\PsFigCap{55}{octahedron}%
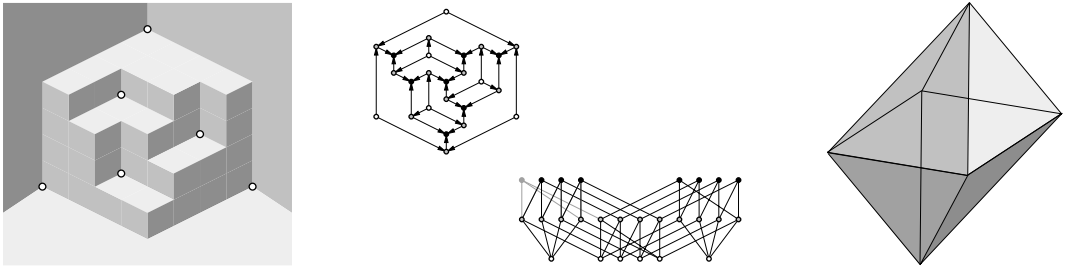
{An orthogonal surfaces, two diagrams of its complex
(save the {\bf 0} element) and the corresponding 3-polytope.}

Figure~\ref{fig:octahedron} shows an example in dimension 3. 
A proof of the theorem is given in~\cite{bps-mr-98}. 
The dimension 3 case of Scarf's theorem was independently obtained by
Schnyder~\cite{s-pgpd-89}. 

It is known that not all simplicial $d$-polytopes have a corresponding
Scarf complex. For example there are neighborly simplicial
4-polytopes, they have the complete graph as skeleton.
From bounds for the dimension of complete graphs it can be concluded
that for $n\geq 13$ these 4-polytopes are not realizable by an
orthogonal surface in
$\RR^4$.

A more general criterion was developed by Agnarsson, Felsner and
Trotter~\cite{aft-mnegb-99}. They show that the number of edges of a
graph of dimension 4 can be at most $\frac{3}{8}n^2 + o(n^2)$.

It would be very interesting to know more criteria which have to be
satisfied by simplicial $d$-polytopes which are Scarf, i.e., can
be realized by an orthogonal surface.

Scarf's Theorem is the tool for the proof of Theorem~\ref{thm:w-span}.
Let $X\subset \RR^2$ be a generic point set with $|X|=n$ and suppose
that all coordinates of points in $X$ are strictly between $0$
and~$M$. We map a point $x=(x_1,x_2)$ in $X$ to the vector $v^x =
(x_1,x_2,M-x_1,M-x_2) \in \RR^4$ and let $V = \{ v^x : x \in X \} \cup
S$, where $S=\{(M,0,0,0),(0,M,0,0),(0,0,M,0),(0,0,0,M)\}$ is the set
of suspension vertices.  The set $V$ is suspended and generic,
therefore, by Scarf's Theorem the Scarf complex corresponds to a
4-polytope. We are interested in the face numbers $F_i$ of the
polytope and consequently also in the face numbers $f_i$ of the
Scarf complex.

\begin{proposition}
The Scarf complex $\Delta_V$ of the above $V$ has face numbers
\Bitem $f_0 = n+4$,
\Bitem $f_1 = \spanned(X) +4n+6$,
\Bitem $f_2 = \span{3}(X) +10n - \exposed(X)$,
\Bitem $f_3 = \span{4}(X) +6n - \exposed(X) -2$.
\end{proposition}

\Proof
There are $n$ vertices $v^x$ and four suspension vertices in $V$.
Therefore, $f_0 = n+4$.

From the characterization ($i$)  we know that $U\subset V$ is a face
of this complex iff there is no vertex in $V$ which is strictly
dominated by $\bigvee\{ v : v\in U \}$. Suppose $U\cap S = \emptyset$,
i.e, $U=\{v^y : y \in Y \}$ for a certain subset $Y$ of $X$,
and let $u=(u_1,u_2,u_3,u_4)$ Since $u_i = \max( v^y_i : y \in Y)$
the definition of $v^y$ implies that $u_i = \max( y_i : y\in Y)$ for
$i=1,2$ and $u_i = \max( M-y_{i-2} : y\in Y) = M - \min(y_{i-2} : y\in
Y)$ for $i=3,4$. It turns out that there is a strictly dominated
vertex $v^x < u$ iff $x$ is contained in the rectangle $R[Y]$ with first
coordinate between $\max( y_1 : y\in Y)$ and $\min(y_1 : y\in Y)$
and second coordinate between 
$\max( y_2 : y\in Y)$ and $\min(y_2 : y\in Y)$, see Figure~\ref{fig:nonempty}.

\PsFigCap{55}{nonempty}%
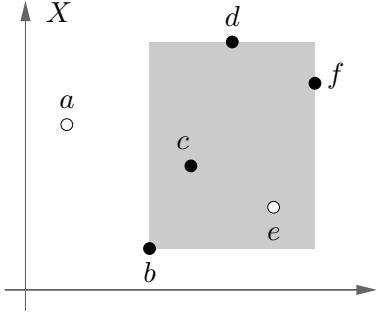
{The rectangle corresponding to $Y=\{b,c,d,f\}$ contains $c$ and $e$ in the 
interior, therefore, $\{v^y : y \in Y\} \not\in \Delta_V$.}

This correspondence explains the occurrences of $\spanned(X)$, $\span{3}(X)$ and
$\span{4}(X)$ in the equations for $f_1$, $f_2$ and $f_3$ given in the proposition.
The remaining terms in these expressions are needed for the
contributions of the suspension vertices. In the case of $f_1$ this is
rather easy. Each pair $(v^x,s)$ with $x\in X$ and $s\in S$ is a
1-face, for example if $s$ is the suspension for coordinate four, then
$v^x\join s = (x_1,x_2,M-x_1,M)$ and $v^y < v^x\join s$ would
require $y_1 < x_1$ and $M-y_1 < M-x_1$ which is impossible.
In addition to these $4n$ 1-faces there are 6 which connect pairs
of suspension vertices. Therefore, $f_1 = \spanned(X) +4n+6$.

We now turn to the 4-element subsets counted by $f_3$, as already
shown there are $\span{4}(X)$ such sets which contain no suspension
vertex.  To count those which contain suspension vertices it is
convenient to adapt the plane visualization with empty rectangles.
Recall that the rectangle $R[Y]$ has its right boundary at $u_1 =
\max( y_1 : y\in Y )$ and this is just the 1st coordinate of
$\bigvee\{ v^y : y\in Y \}$. Similarly the top, left and bottom of
$R[Y]$ are at $u_2$, $u_3$ and $u_4$.  Now suppose that the first
coordinate of $\bigvee\{ v : v\in U \}$ is $M$, this can be
represented by allowing $R$ an unbounded right side. Again $R$
corresponds to a face of the Scarf complex iff it contains no points
of $X$ in the interior.

Now we can count these generalized rectangles.  To begin with we find
$4n$ generalized rectangles by starting at a point $x\in X$ and
selecting one of the four directional rays starting from $x$. This ray
is then widened to both sides until it hits a point from $X$ (see
Figure~\ref{fig:special3faces} for some widened rays) . If the widened
ray doesn't hit a point there is another unbounded side.  Note that
this procedure counts all empty generalized rectangles with just one
and the four with exactly three unbounded sides, see
Figure~\ref{fig:special3faces}.  Those with exactly two unbounded
sides which are generated are generated twice. There are as many of
this kind as there are edges on the ortho-hull, this number in turn
equals the number of points on this hull which is just the number
$\exposed(X)$. Still uncounted are those rectangles which have exactly
two unbounded sides opposite to each other, the figure shows an
example. There are $2(n-1)$ rectangles of that type. Collecting the
numbers we find $f_3 = \span{4}(X) + 4n - \exposed(X) + 2n - 2$.

For the 3-element subsets counted by $f_2$, we know that
$\span{3}(X)$ such sets contain no suspension
vertex. Each triple of suspension vertices contributes to $f_2$.
For those triples containing one or two points from $X$ we again consider generalized
empty triangles. For each point $x$ and each of the four directional
rays, the ray can be widened to either side until it hits a point from $X$, if
it doesn't hit a point there is another unbounded side. This gives a
total of $8n$ but we have generated rectangles with two unbounded
sides twice. There is one such generalized rectangle associated to
each exposed point, the four points which have a extreme coordinate
actually have two empty quadrants. Thus the corrected contribution is
$8n - \exposed -4$. Still uncounted are those rectangles which have exactly
two unbounded sides opposite to each other. These are $2n$ since each point of $X$
belongs to two of them . Summing up we obtain
$f_2 = \span{3}(X) + 4 + 8n  - \exposed(X) - 4 + 2n$.
\qed

\PsFigCap{40}{special3faces}%
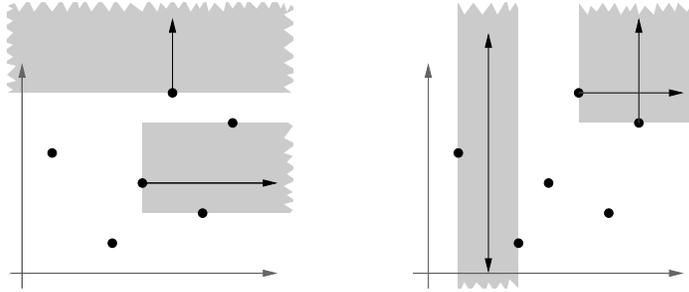
{Generalized empty rectangles, left types with one or three unbounded
  sides, right the two types with two unbounded sides. }

Given this proposition the proof of Theorem~\ref{thm:w-span} is 
easy. 
\medskip

\Proof (Theorem~\ref{thm:w-span})
Let $F_i$ be the face numbers of the 4-polytope corresponding to the
Scarf complex $S_V$. This polytope has one facet which is not
represented in the Scarf complex, this is the facet of the four
suspension vertices, i.e., $F_3 = f_3 +1$ and $F_i = f_i$ for
$i=0,1,2$.

The Euler-Poincar\'e formula $1 - F_0 + F_1 - F_2 + F_3 - 1 = 0$
together with the double counting identity for simplicial 4-polytopes
$2F_2=4F_3$ yields $F_3 = F_1 - F_0 = f_1 - f_0 = \spanned + 3n + 2$.
On the other hand $F_3 = f_3 +1 = \weaksp(X) + 6n - \exposed(X) -1$.
Combining the two expressions completes the proof of the first
relation.

For the second relation of the theorem use Euler-Poincar\'e to obtain
$F_2 = 2\,(F_1 - F_0) = 2\,\spanned + 6n + 4$ and combine it with the
expression for $f_2$ from the proposition.
\qed



\bibliography{os-bib}
\bibliographystyle{/homes/combi/felsner/FU-DATA/TexStuff/my-siam}


\newpage

\section{Appendix: Proof of Theorem~\ref{thm:span}}\label{sec:proof1}

Let $X$ be a set of $n$ points in $\RR^2$ with a maximal number of
empty spanned rectangles. Let $G_r(X)$ be the rectangle graph on $X$.
We have to prove that $|E_r(X)| \leq \lfloor\frac{n^2}{4} +n -2\rfloor$.

Let $m$ and $M$ be the points with minimal and maximal first
coordinate. If the number of edges incident to these two points is at
most $(n-2) + 3$, then we can remove them and do with
induction: $|E_r(X)| \leq |E_r(X-\{m,M\})| + (n-2) + 3 \leq
\big\lfloor\frac{(n-2)^2}{4} +(n-2) -2\big\rfloor + (n-2) + 3 =
\lfloor\frac{n^2}{4} +n -2\rfloor$.

Now suppose that $m$ and $M$ are incident to at least $(n-2)+4$ edges.
This implies that the two points have at least 3 common neighbors.
The following lemma is useful to localize the common neighbors.

\begin{lemma}\label{lem:streifen}
Let $(x,y)$ be an edge and let $x=(x_1,x_2)$ and $y=(y_1,y_2)$,
then there is at most one common neighbor of $x$ and $y$ in
each of the four region $A^\uparrow(x,y)$, $A^\rightarrow(x,y)$,
$A^\downarrow(x,y)$ and $A^\leftarrow(x,y)$ shown in Figure~\ref{fig:streifen}.
\end{lemma}
\PsFigCap{55}{streifen}%
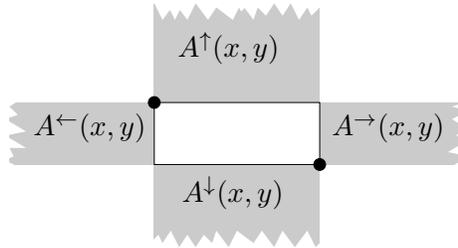
{Areas in which $x$ and $y$ can only have one common neighbor.}
\Proof
We prove the statement for the region $A^\uparrow(x,y)$.
Let $z=(z_1,z_2) \in A^\uparrow(x,y)$ i.e.,
$\min(x_1,y_1) < z_1 < \max(x_1,y_1)$ and $z_2 > \max(x_1,y_1)$.
Suppose $z$ is a neighbor of $x$ and $y$ and consider
$z'$ with $z_2' > z_2$. Either
$z \in  R[x,z']$ or $z\in R[y,z']$, i.e., one of rectangles is
non-empty. Therefore, $z'$ can not be a common
neighbor of $x$ and $y$.\qed

Since $A^\uparrow(m,M)$ and $A^\downarrow(m,M)$
contain at most one common neighbor of $m$ and $M$ and 
by the choice of the points $A^\rightarrow(m,M) = \emptyset$
and $A^\leftarrow(m,M)=\emptyset$ it follows that $R[m,M]$ can not be
empty. This implies that $m,M$ is not an edge and there must be at least
two common neighbors of $m$ and $M$ in $R[m,M]$.
Let $a$ and $b$ be two such points. Without loss of generality we
may assume $M_2 > m_2$ and  $a_1 < b_1$.

If the number of edges incident to the four points $\{m,M,a,b\}$ is at
most $2(n-4) + 8$, then we can remove these points and complete with
induction. The four points already induce at least 4 edges, if
$R_0=R[a,b]$ is empty they even induce 5.  Figure~\ref{fig:four-points}
shall help us analyze where we can expect points which are common
neighbors to more than two from the set $\{m,M,a,b\}$. Given a
3-subset of $\{m,M,a,b\}$ the region where common neighbors of the
three points can live are restricted by the emptiness of the spanned
rectangles.

\PsFigCap{70}{four-points}%
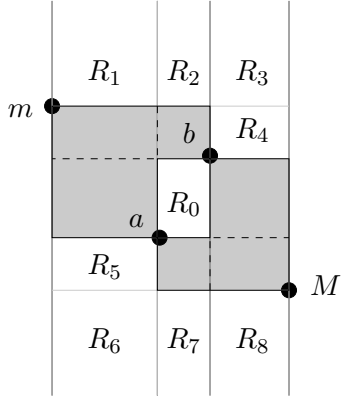
{Regions defined by position relative to points $m,M,a,b$.}

\Claim 1. $\{m,a,b\}$ can have at most two common neighbors,
one in $R_0$ and one in $R_1\cup R_2$.
\Proof Regions $R_3$ and $R_4$ are obstructed for $a$ by $b$.
Regions $R_5$ and $R_6$ are obstructed for $b$ by $a$ and
Regions $R_7$ and $R_8$ are obstructed for $m$ by $a$.
Since $R_0 = A^\rightarrow(m,a)$ and $R_1\cup R_2 = A^\uparrow(m,b)$
Lemma~\ref{lem:streifen} implies the statement.\qedclaim

The proofs of the following claims are similar and omitted.

\Claim 2. $\{m,a,M\}$ can have at most two common neighbors,
one in $R_0$ and one in $R_5\cup R_6$.

\Claim 3. $\{m,b,M\}$ can have at most two common neighbors,
one in $R_0$ and one in $R_3\cup R_4$.

\Claim 4. $\{a,b,M\}$ can have at most two common neighbors,
one in $R_0$ and one in $R_7\cup R_8$.
\smallskip

\ni
This makes a total contribution of 8 edges from the points of high
degree. Together with the 4 or 5 edges induced by the four selected
points this is 12, way too much.
The following claim reduces the contribution of $R_0$ to 2 and, hence,
the contribution from high degree points to~6.

\Claim A. if $R_0$ contains $r_0$ points, then these
points can share at most $2r_0 + 2$ edges with $\{m,a,b,M\}$.
\medskip

Let $z\in R_0$ be a point with more than two of these edges,
suppose $a,m,b$ are neighbors of $z$. It follows that $z$ is
the only point which spans an empty rectangle with the middle
of the three, in our case this is $m$. The only remaining triple is
$a,M,b$ but again, this triple can only contribute 1.\qedclaim

If the contribution from high degree points is 6, then
there are points in $R_0$ as well as in $R_3\cup R_4$ and in
$R_5\cup R_6$. When we remove $a$ and $b$ these points will
span at least two new empty rectangles which are not present
in the original. These two additional edges may be subtracted
from the inductive estimate. This allows us to get around with
a total of $2(n-4) + 4 +6$ edges incident to points in $\{m,a,b,M\}$.

The situation where the contribution from the high degree points is 5
is similar. We can, as before,
take advantage of empty rectangles that appear when removing
 $a$ and $b$.

\PsFigCap{30}{best}%
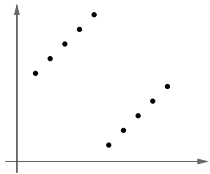
{A point set with 10 points maximizing $\spanned$.}

The example shown in Figure~\ref{fig:best} indicates a construction
which yields a set of points with exactly $\lfloor\frac{n^2}{4} + n -
2\rfloor$ edges for every $n\geq 2$.
\qed

\section{Appendix: Proof of Theorem~\ref{thm:singlearrow}}\label{sec:proof2}

As in the previous proof we use induction. We remove some vertices
such that the number of edges incident to them is small enough.
The next lemma quantifies the term `small enough'.

\begin{lemma}\label{lem:induct}
Let $\cal G$ be a class of graphs such that for every member $G$ of
$\cal G$ with $n$ vertices it is possible to remove a set of 
$2k$ vertices which is incident to at most $k(n-2k) + 2ks$ edges,
then the number of edges of graphs in $\cal G$ is
at most $\frac{1}{4}n^2 + sn$.
\end{lemma}

\Proof
Let $e_n$ be the maximal number of edges of a graph with $n$
vertices in $\cal G$. By induction 
$e_n \leq k(n-2k) + 2ks + e_{n-2k} \leq  
k(n-2k) + 2ks + \frac{1}{4}(n-2k)^2 + s(n-2k)$.
This yields $e_n \leq \frac{1}{4}n^2 + sn  - k^2$.\qedclaim

Let $G=(V,E)$ be a graph with dimension at most $\between{3}{4}$
and let $\pi_1,\pi_2,\pi_3,\pi_4=\ol{\pi_3}$ be a realizer of $G$.
For a vertex $v$ let $v_i$ be the position of $v$ in $\pi_i$.
Condition $(*)$ is that for every edge $u,v$ and $w\neq u,v$
there is an $i$ with $w_i > \max(u_i,v_i)$.
Writing $[s,t]$ for the interval of integers between $s$ and $t$
we can restate the edge condition as follows:
$u,v$ is admissible for an edge if every $w \neq u,v$
fulfills $w_1 > \max(u_1,v_1)$ or $w_2 > \max(u_2,v_2)$ or
$w_3 \not\in [u_3,v_3]$.

For points $p,q\in \RR^2$ let $\RN[p,q] = R[0,p,q]$ be the rectangle 
spanned by the two points together with the origin.  
The first two `coordinates' of a vertex $v\in V$ give
a point in the plane $P_{1,2}$ (we introduce no new notation and simply call
the point $v$). Since we are interested in maximizing the number of
edges we can assume that every pair $u,v\in V$ which is admissible for
an edge by condition $(*)$ is actually in $E$. 
Here is a geometric description of the edges of~$G$:

{\Item($\star$) A pair $u,v\in V$ is an edge of $G$ iff every $w \in \RN[u,v]$
satisfies $w_3 \not\in [u_3,v_3]$.
}\par

\ni
Given a vertex $v$ we classify the neighbors of $v$ according to the
quadrant of $v$ containing them in the plane $P_{1,2}$. Formally:
\begin{eqnarray*}
N^\ne(v) & = &
\{ u \in N(v) : u_1 > v_1 \text{ and } u_2 > v_2 \}\\
N^\nw(v)  & = &
\{ u \in N(v) : u_1 < v_1 \text{ and } u_2 > v_2 \}\\
N^\sw(v)  & = &
\{ u \in N(v) : u_1 < v_1 \text{ and } u_2 < v_2 \}\\
N^\se(v)  & = &
\{ u \in N(v) : u_1 > v_1 \text{ and } u_2 < v_2 \}
\end{eqnarray*}

\Fact 1. A vertex $v$ has at most $2$ neighbors in $N^\sw(v)$ 
\medskip

\Proof The two candidates for $u\in N^\sw(v)$ are the vertices with
$u_3$ maximal subject to the condition $u_3 < v_3$ and with
$u_3$ minimal subject to $u_3 > v_3$. These vertices obstruct
the $\star$ property for all the other vertices in this quadrant 
of $v$.
\qedclaim
 
To simplify the analysis we will disregard all the edges $\{u,v\}$ with
$u \in N^\sw(v)$. Thereby we loose at most
$2n$ edges which is insignificant since the bound stated in the
theorem allows imprecision of order $O(n)$. Let $G'=(V,E')$ be the
remaining graph.

All neighbors of a vertex $v$ in $G'$ are either north-west or south-east
of $v$. It is useful to consider a finer partition of the neighborhood $N(v)$
of $v$:  
\begin{eqnarray*}
N^\nw_+(v) = \{ u \in N^\nw(v) :  u_3 > v_3 \} & & 
N^\se_+(v) = \{ u \in N^\se(v) :  u_3 > v_3 \}\\[3pt]
N^\nw_-(v) = \{ u \in N^\nw(v) :  u_3 < v_3 \} & &
N^\se_-(v) = \{ u \in N^\se(v) :  u_3 < v_3 \}
\end{eqnarray*}

\ni
Each of these sets has the nice property that among the members of the set two
of the coordinates are bound together:

\Fact 2.\qquad\raise\baselineskip\hbox{\parbox[t]{.6\hsize}{
\Item(a) $w,w'\in N^\nw_+(v)$, then $w_2 < w'_2 \iff w_3 > w'_3$,
\Item(b) $w,w'\in N^\se_+(v)$, then $w_1 < w'_1 \iff w_3 > w'_3$,
\Item(c) $w,w'\in N^\nw_-(v)$, then $w_2 < w'_2 \iff w_3 < w'_3$,
\Item(d) $w,w'\in N^\se_-(v)$, then $w_1 < w'_1 \iff w_3 < w'_3$.
}}\medskip

\Proof
We only prove (a), the other statements are obtained by permutations in coordinates and signs.
Let $w,w'\in N^\nw_+(v)$ with $w_2 < w'_2$ this implies that $w \in \RN[w',v]$,
hence, $w_3 \not\in[v_3,w'_3]$. This together with $w\in N_+(v)$, i.e., $w_3 > v_3$,
implies $w_3 > w'_3$. Conversely, $w,w'\in N^\nw_+(v)$ and $w_3 > w'_3$
implies $v_3 < w'_3 < w_3$. Since $w\in N(v)$ this requires $w'\not\in \RN[w,v]$
which forces $w_2 < w'_2$.\qedclaim
 
The common neighbors of a pair $u,v$ of vertices are located north-west or south-east
of both or they lie in the rectangle spanned by $u$ and $v$:
\begin{center}
$N^\nw(u,v) =  N^\nw(u) \cap N^\nw(v) \qquad
N^\se(u,v) =  N^\se(u) \cap N^\se(v)$ \\[4pt]
$N^\rec(u,v)   =  \{ w \in N(u)\cap N(v)  : w \in R[u,v] \}$
\end{center}

In the following analysis we say that a set $A$ is \term{essentially contained} in $B$,
denoted by $A \subseteq^* B$ if this is true with at most four exceptional elements,
i.e., $ (A\setminus A') \subseteq B$ for some $A'$ with $|A'|\leq 4$.

\Fact 3.\qquad\raise\baselineskip\hbox{\parbox[t]{.6\hsize}{
\Item(a) If $u_1 < v_1$ and $u_3 < v_3$, then $N^\nw(u,v) \subseteq^* N^\nw_+(v)$.
\Item(b) If $u_1 < v_1$ and $u_3 > v_3$, then $N^\nw(u,v) \subseteq^* N^\nw_-(v)$.
\Item(c) If $u_2 > v_2$ and $u_3 > v_3$, then $N^\se(u,v) \subseteq^* N^\se_+(u)$.
\Item(d) If $u_2 > v_2$ and $u_3 < v_3$, then $N^\se(u,v) \subseteq^* N^\se_-(u)$.
}}\medskip

\Proof 
We only prove (a), essentially the same argument, with appropriate
permutations of coordinates shows (b), (c) and (d).

If $w \in N^\nw(u,v) = N^\nw(u)\cap N^\nw(v) $ then $u \in \RN[w,v]$.
Therefore, an edge $w,v$ requires $w_3 > u_3$. There can only be one
vertex $w$ in $N^\nw(u,v)$ with $u_3 < w_3< v_3$ since $w$ obstructs
either the edge with $u$ or the edge with $v$ for every $w'$ with
larger 2nd coordinate. Hence, all but at most one element from
$N^\nw(u,v)$ have 3rd coordinate larger than $v$, i.e., are members of
$N^\nw_+(v)$. \qedclaim

Let $m$ and~$M$ be the elements of $V$ with maximal and minimal 3rd coordinate.
If $m$ and~$M$ are incident to at most $(n-2) + 6$ edges in $E'$, then
we can apply Lemma~\ref{lem:induct} with $k=1$ and $s=3$.

From Fact~3 it follows that the essential portion of common
neighbors of $m$ and $M$ must be in $N^\rec(m,M)$. In particular, if there
are more than two common neighbors of $m$ and $M$ in $G'$, then they
are the south-east and the north-west corners of their rectangle $R[m,M]$.
By symmetry in the first two coordinates we may assume
that $m$ is east, i.e., $m_1 > M_1$ and $m_2 < M_2$.

Combining part (a) and (d) of Fact~2 it follows that for $w,w'\in
N^\rec(m,M)$ the coordinates are bound by 
\begin{equation}\label{eq:N(mM)}
w_1 < w'_1 \iff w_2 > w'_2
\iff w_3 < w'_3.
\end{equation}
Let $a$ and $A$ be the elements of $N^\rec(m,M)$
with minimal and maximal 3rd coordinate. Figure~\ref{fig:MaAm} shows
the relative positions of $M,a,A,m$ in the planes $P_{1,2}$ and
$P_{1,3}$. Now consider the edges incident to these four vertices in
$E'$. If their number is at most $2(n-4) + 12$, then we apply
Lemma~\ref{lem:induct} ($k=2$, $s=3$).

The next goal is to identify points which are neighbors to at least three vertices
from $M,a,A,m$. The common neighbors of $m$ and $M$ have been localized before.
From the way their coordinates are coupled (see~\ref{eq:N(mM)}) it follows that 
both of $a$ and $A$ can have at most one common neighbor with $m$ and $M$.
Therefore, we have to look for common neighbors of $a,A,m$ or $a,A,M$.
Appealing to symmetry we concentrate on the first of these cases.

Fact~3(b) implies that $N^\nw(A,m) \subseteq^* N^\nw_-(m)$ 
but $N^\nw_-(m)= \emptyset$ by the choice of $m$.
Fact~3(c) implies that $N^\se(A,m) \subseteq^* N^\se_+(A)$ 
but $N^\se(a,A) \subseteq^* N^\se_-(a)$ by 3(d). That is, all
but at most one of the common neighbors of $a,A$ and $m$ are in
$N^\rec(A,m)$. These elements are in $N^\se(a,A) \subseteq^*
N^\se_-(a)$ and in $N^\nw_+(m)$. 
From part (a) and (d) of Fact~2 it follows that for 
two common neighbors $w,w'$ of $a,A,m$ the coordinates are bound by 
$w_1 < w'_1 \iff w_2 > w'_2 \iff w_3 < w'_3$, i.e.,
as in~\ref{eq:N(mM)}.

Let $b$ and $B$ be the common neighbors with minimal and maximal 3rd
coordinate (see Figure~\ref{fig:MaAm}). Now consider the edges
incident to the vertices $M,a,A,b,B,m$ in $E'$. If their number is at
most $3(n-6) + 18$, then we apply Lemma~\ref{lem:induct} ($k=3$,
$s=3$).
\PsFigCap{58}{MaAm}%
{The relative positions of $M,a,A,b,B,m$ in the planes $P_{1,2}$ and
$P_{1,3}$.}

If there are too many incident edges, then there are `many' vertices 
which are adjacent to at least four out of the six vertices.
As before $m$ and $M$ can not both belong to such a four set.
Consider  four sets containing $a$ and $M$, an argument as before
shows that all
but at most one of the common neighbors of $a$ and $M$ are in
$N^\rec(M,a)$ and they have 3rd coordinate larger than $a_3$.
Therefore, $a$ obstructs edges from these vertices to $b,B$
and $m$. Similarly, the pairs $A,M$, $b,m$ and $B,m$ can not
be contained in a four set with many common neighbors.
The only remaining candidate for such a four set is $a,A,b,B$.
These elements may indeed have a large common neighborhood.
These neighbors must be contained in 
$N^\rec(A,b)$ because they are in $N^\se(a,A) \subseteq^*
N^\se_-(a)$ and in $N^\nw(b,B) \subseteq^* N^\nw_+(B)$. 
Hence, we again know that for any two  $w,w'$ of them
the coordinates are again bound as in ~\ref{eq:N(mM)}:
$ w_1 < w'_1 \iff w_2 > w'_2 \iff w_3 < w'_3$.

Let $c$ and $C$ be the common neighbors of $a,A,b,B$ with minimal and
maximal 3rd coordinate. Now consider the edges incident to 
the vertices $M,a,A,c,C,b,B,m$ in
$E'$. If their number is at most $4(n-8) + 24$, then we are done by
Lemma~\ref{lem:induct} ($k=4$, $s=3$). 

If there were more edges incident to these eight vertices, then there
would have to be `many' vertices 
which are adjacent to at least five out of the eight vertices.
Along the lines of the previous analysis, however, it can be shown
that there is no five element subset of vertices from
$M,a,A,c,C,b,B,m$ which has more than four common neighbors.
The set itself induces only 14 edges, so that $4(n-8) + 24$
is truly an upper bound for the number of edges incident to
vertices from $M,a,A,c,C,b,B,m$.
\qed

\end{document}